\documentclass[12pt]{amsart}

\usepackage{common}
\usepackage{mathrsfs}

\renewcommand{\fS}{\m{\frak S}}
\newcommand{\U}{\m{\frak U}}

\title{Generic separable metric structures}

\author{Alexander Usvyatsov}

\address{Alexander Usvyatsov\\
  University of California -- Los Angeles\\
  Mathematics Department\\
  Box 951555\\
  Los Angeles, CA 90095-1555\\
  USA}

\urladdr{http://www.math.ucla.edu/\textasciitilde alexus}

\date{\today}

\begin{document}

\begin{abstract}
  We compare three notions of genericity of separable metric structures. Our
  analysis provides a general model theoretic technique of showing that structures are
  generic in descriptive set theoretic (topological) sense and in measure theoretic sense.
  In particular, it gives a new perspective on Vershik's theorems on
  genericity and randomness of Urysohn's space among separable metric spaces.
\end{abstract}

\maketitle

\section{Introduction}
There are several ways to define the notion of a ``generic'' metric structure. In this
article we compare the model theoretic and two topological approaches to
this question.

This work was motivated by Anatoly Vershik's results on genericity
and randomness of the Urysohn space among separable metric spaces,
Theorems 1 and 2 in  \cite{Ver02}. Vershik considers the collection
of all separable metric spaces as a topological space, let us call
it $\fS$. Some elements of $\fS$ are (isometric to) the Urysohn
space. Vershik shows that this set is $G_\delta$ dense in $\fS$,
which leads to the conclusion that the Urysohn space is in a sense
``a generic'' separable metric space. Then he shows that for any
``reasonable'' probability measure on $\fS$, the collection of
metric spaces isometric to the Urysohn space is of measure $1$. This
leads to the conclusion that the  Urysohn space is in a sense ``the
random'' metric space.

In his talk at the workshop on the Urysohn space at Ben-Gurion University
(May 2006), Vershik said that his results had
been motivated by model theoretic properties of the (countable) random graph, and that the
theorems in \cite{Ver02} are in some sense the analogues of the appropriate facts in classical model
theory, although the context is different: instead of countable structures one deals
with topological spaces of cardinality the continuum. In this paper we aim to show that the
analogy goes much farther.

Indeed, \emph{countable discrete} structures are replaced in this context
with \emph{separable metric} spaces, so classical model theory is not the appropriate general
framework. We would like to convince the reader that there exists a natural generalization of
\emph{discrete} first order logic to the \emph{continuous} context, in which
Vershik's results are the true analogues of classical facts, the Urysohn space is an
analogue of the random graph, and discrete countable models are no more than a
particular case of separable continuous structures. So from our point of
view, properties of the Urysohn space discovered by Vershik are much more than results
inspired by certain similarities between this structure and the random graph; in a sense,
both of these are particular cases of the same model theoretic phenomenon, which we intend
to describe here.

Continuous first order logic, recently introduced by Ita\" i Ben-Yaacov
and the author in \cite{BU}, allows one to study classes of metric spaces
(maybe equipped with continuous extra-structure, e.g. a collection of uniformly continuous
functions from the spaces to $\setR$) from model theoretic point of view. Once working in this context,
many results in classical model theory generalize to analytic structures. This paper is
devoted to the connection between model theory and descriptive set theory, which is
very well-developed in the classical context, i.e. studying Polish spaces of countable
structures for a given countable signature, countable models of a countable universal theory,
etc. We will refer the reader to the excellent expository paper by Greg Hjorth, \cite{Hj04}.

In addition to generalizing Vershik's theorem to a broad collection
of classes of metric structures, our work generalizes a few basic
concepts and results from classical first order model theory to the
continuous context and pushes out the boundaries of possible
applications of continuous logic. So although we intentionally try
to make the article accessible to non-logicians, it could also be of
interest to model theorists.

Working in the context of continuous first model theory, we adapt some basic facts and
techniques from \cite{Hj04} and show how one defines a Polish topology on the space of e.g.
all separable models of a certain universal continuous theory. Having done that, we discuss
three different notions of genericity of a structure. One is model theoretic, genericity of a model
of a universal theory among its peers). The other two are topological,
genericity of a structure as an element of the appropriate Polish space in two different ways:
in the sense of Baire category theory and in the sense of measure theory. Let us
state things more precisely.


Let $K$ be a ``reasonable'' class of separable metric structures.
In our context $K$ will normally be the
class of all separable models of a certain universal continuous first order theory.
From the model theoretic point of view, a generic structure in $K$ is a structure
in which ``anything that can happen'' in $K$ happens. Such structures are called ``existentially
closed'' for $K$. We will give precise definitions later.

On the other hand, one can consider $K$ as a Polish space
(i.e. there is a natural topology on $K$ with respect to which $K$ is a complete separable
metric space).
One can call a structure ``generic''
for $K$ in topological sense if its isomorphism class
is a ``big''
subset of $K$. One natural notion of ``bigness'' in this context is $G_\delta$ dense.
Another one comes from measure theory: one can consider natural measures on the space $K$
and ask what are the sets of measure $1$.

In this article we have several primary goals:
\begin{enumerate}
\item Introduce the general model theoretic framework and the relevant notion
  of genericity.
\item Construct the Polish space of separable metric structures.
\item Connect the notions of genericity. More precisely, we explain
  how a model theoretic notion of genericity gives rise to $G_\delta$ dense
  sets in the appropriate Polish space $\fS$ and
  sets of measure $1$ with respect to any ``reasonable''
  probability measure on $\fS$. In particular, this provides a
  powerful general technique for showing that certain structures are
  topologically generic and random (as it allows us to use well-developed model
  theoretic tools for this purpose).
\item
  Conclude with some examples. In particular, we discuss
  model theory of Urysohn space and
  show that our results generalize Vershik's theorems on its ``topological'' genericity.
\end{enumerate}

\subsection*{Acknowledgements} The author thanks the anonymous referee for very helpful
comments and suggestions.

\section{Preliminaries and basics}
\subsection {Continuous logic}

Continuous first order logic was introduced in \cite{BU} and
developed further by Ita\"{i} Ben-Yaacov, Alexander Berenstein,
C. Ward Henson and the author. We refer the reader to \cite{BBHU} for
a detailed exposition. We will now try to summarize some important basic notions,
facts and notations.

Just as in classical predicate logic,
one starts with a fixed \emph{signature (vocabulary)} $\tau$. In this paper, $\tau$ will be
countable. A signature (vocabulary) is a collection of function symbols and predicate symbols
as well as continuity moduli for all these symbols. There is a distinguished predicate
symbol $d(x,y)$, which will correspond to the metric.

Given a vocabulary $\tau$, one constructs the continuous \emph{language} $L$ which corresponds
to it, which consists of continuous first order $\tau$-\emph{formulae}.
As in classical first order logic, formulae are constructed by induction using \emph{connectives}
and \emph{quantifiers}. Any countable collection of continuous functions from $[0,1]^k$ to $[0,1]$
(for any $k$) which is dense in the set of all such continuous functions
can be taken as the set of connectives. We will assume that
the following functions are among our connectives:

\begin{enumerate}
\item
  The \emph{constant function} $q$ for every $q\in [0,1]\cap\setQ$
\item
  \emph{pointwise minimum} ($[0,1]^2 \to [0,1]$)
\item
  \emph{pointwise maximum} ($[0,1]^2 \to [0,1]$)
\item
  \emph{Multiplication by} $q$, $[x \mapsto x\cdot q]$, for every $q\in [0,1]\cap\setQ$ ($[0,1]\to [0,1]$)
\item
  \emph{negation}, $[x \mapsto 1-x]$, ($[0,1]\to [0,1]$)
\item
  \emph{dotminus} or \emph{implication}: Truncated (at 0) minus $[(x,y) \mapsto x \dotminus y]$  ($[0,1]^2 \to [0,1]$)
\item
  Truncated (at 1) \emph{plus} $[(x,y) \mapsto x + y]$  ($[0,1]^2 \to [0,1]$)
\item
  $(x,y) \mapsto |x-y|$  ($[0,1]^2 \to [0,1]$)
\end{enumerate}

Of course, some of the functions above can be defined using the others, but we are not looking for
``minimal'' systems of connectives here.

The continuous \emph{quantifiers}
are $\inf_x$ and $\sup_x$. As in classical first order logic, we only allow quantification
over elements.

So the following are examples of formulae:
\begin{itemize}
\item $d(x,y)$
\item $d(x,y)\dotminus d(y,x)$
\item $\inf_{x,y}d(x,y)$
\item $\sup_{x,y,z}(d(x,z)\dotminus(d(x,y)+d(y,z)))$
\end{itemize}

As usual, formulae with no ``free variables'' (i.e. each variable is in a scope of
one of the quantifiers) are called \emph{sentences}. The first two formulae above are not
sentences, while the last two are.

An $L$-\emph{pre-structure} is a set $M$ equipped with \emph{interpretations} for all $\tau$-symbols
such that $d$ is interpreted as a pseudometric, each predicate symbol is interpreted as a function
from (some power of) $M$ to $[0,1]$, each function symbol is interpreted as a function from
(some power of) $M$ to $M$, and all of them respect their continuity moduli with respect to $d$. In
other words:

\begin{itemize}
\item
  $d^M\colon M^2 \to [0,1]$ is a pseudometric
\item
  For every $n$-ary predicate symbol $P$, we have $P^M \colon M^n \to [0,1]$
  uniformly continuous
  with respect to $d$ (respecting the continuity modulus of $P$ dictated by $\tau$)
\item
  For every $n$-ary function symbol $f$, we have $f^M \colon M^n \to M$
  uniformly continuous
  with respect to $d$ (respecting the continuity modulus of $f$ dictated by $\tau$)
\end{itemize}

A \emph{structure} is a pre-structure in which $d$ is a complete metric.

See \cite{BU} or \cite{BBHU} for more details (on e.g. continuity moduli). Formal definitions
of these notions are not important for us here; but it is crucial that the interpretation of
each predicate
symbol and of each function symbol
is uniformly continuous, and uniformly so in all $L$-structures (this is what we
need the continuity moduli for). Uniform continuity allows us to take ultraproducts of $L$-structures
and obtain e.g. compactness of first order continuous logic.

Note that given a structure $M$, one can easily define (by induction)
the $M$-value of $\phh$ for any sentence $\ph$. We will denote this value
by $\ph^M$ (it is a real number in the interval $[0,1]$).

Note also that there is no particular importance for the interval $[0,1]$, but
every predicate symbol must have bounded range (again, so that
ultraproducts will work), and by rescaling we may assume it is in fact always $[0,1]$.

A \emph{condition} is a statement concerning the value of a sentence $\ph$. For example,
$\ph\le\eps$, $\ph=0$, $\ph < \eps$ are conditions (where $\eps \in [0,1]$).
We will call conditions of the form $\ph\le\eps$, $\ph=0$, etc
\emph{closed conditions} and those of the form $\ph < \eps$, etc \emph{open conditions}.

Note that as we will mostly work with conditions of the form $\ph \le \eps$ and $\ph < \eps$,
the continuous quantifiers $\inf_{x}$ and $\sup_{x}$ can be viewed as analogues of the existential and
the universal quantifiers respectively.

It is clear what it means for a structure $M$ to \emph{satisfy} a condition $\al$, and
we write $M \models \al$. If
every structure which satisfies \all also satisfies \bee, we say that \bee
\emph{follows} from \all and write $\al \models \be$. If $\Lambda$ is a collection of
conditions and $M$ is a structure, we say that $M$ is a model for (of) $\Lambda$ if
$M$ satisfies every condition in $\Lambda$, and write $M\models \Lambda$.

A \emph{theory} $T$ is a collection of \emph{closed} conditions which is consistent (i.e. there is a
structure $M$ which satisfies all the conditions in $T$, $M \models T$). We will always assume that
theories are closed under entailment, i.e. if $\al\in T$ and $\al \models \be$, then
$\be \in T$. We denote by $\Mod(T)$ the class of all models of $T$.


We encourage the reader to have a look at examples of continuous languages and theories presented
in \cite{BU} and \cite{BBHU}.

We shall not discuss ultraproduct constructions in this paper. Again, curious readers are
referred to \cite{BU} or \cite{BBHU}. An important consequence is the \emph{Compactness
Theorem} for continuous logic, which will be useful for us:

\begin{fct}(Compactness Theorem)
    Let $\Lambda$ be a collection of \emph{closed} conditions which is finitely satisfiable
    (i.e. every finite subset of $\Lambda$ has a model). Then $\Lambda$ has a model.
\end{fct}

Let $M \subseteq N$ be $L$-structures. We say that $M$ is an
\emph{elementary submodel} of $N$ ($M\prec N$) if for every
$L$-sentence $\phh$ we have $\ph^M = \ph^N$. We say that a theory
$T$ is \emph{model complete} if for every $M,N \models T$, $M
\subseteq N \then M \prec N$. Most theories are not model complete;
we will discuss this notion more later. $T$ is model complete if
(but not only if) it \emph{eliminates quantifiers}; see more in
\cite{BU} or \cite{BBHU}.

Note that continuous first order logic is a natural generalization of classical first order logic.
Indeed, every classical first order theory can be viewed as a continuous theory in which the metric
is discrete.

\subsection{Polish space of separable continuous structures}

Let $\tau$ be a fixed countable continuous vocabulary. For
simplicity we assume that $\tau$ is relational (i.e. no function symbols). Let
 $L$ be the corresponding (countable) continuous language.
We denote the space of all $L$-continuous separable structures $M$ with a distinguished countable dense subset
$\setN \subseteq M$ by $\fS$. Consider the following topology on $\fS$: basic open sets are
of the form $U_{\ph(\x),\a,\eps} = U_{\ph(\a),\eps} = \set{M \in \fS \colon \ph^M(\a)<\eps}$
where $\ph(\x)$ is a \emph{quantifier free} $L$-formula, $\a \in \setN$,
$\eps \in [0,1]\cup\set{\infty}$.

\begin{prp}
  $\fS$ with the topology above is a Polish space.
\end{prp}
\begin{prf}
  Let \lseq{R}{i}{\om} be an enumeration of $\tau$, $R_0$ being the metric. Let $k_i$ be the arity
  of $R_i$ (so $k_0 = 2$).

  By section 2 of \cite{Hj04}
  the product space
  $$X = [0,1]^{\bigsqcup_i\setN^{k_i}}$$
  is Polish.
  We can view $\fS$ as a subspace of $X$
  via the following embedding $\phi:\fS\to X$:

  $\phi(M) = \lseq{f}{i}{\om}$ such that $f_i$ is precisely $R_i^M$ on
  the dense subset $\setN$ of $M$.

  Note that the fact that $R_0$ is a pseudometric and all the rest of the predicates
  respect the appropriate continuity moduli with respect to it is a collection
  of closed conditions. The fact that $R_0$ is an actual metric can be expressed as a
  collection of open conditions. So $\fS$ can be viewed as a $G_\delta$ subset of
  a Polish space, and therefore, by Lemma 2.2 in \cite{Hj04}, $\fS$ is a Polish space
  itself.
\end{prf}

Let $T$ be an $L$-theory. We denote the space of all elements of $\fS$ which are models of $T$ by
$\fS_T$. So $\fS = \fS_\emptyset$.

\subsection{Universal theories and existentially closed models}

\begin{dfn}
\begin{enumerate}
\item
  We call a theory \emph{universal} if it is (the closure under entailment of) a collection of
  conditions of the form $[\sup_{\x}\ph(\x)=0]$ where $\phh$ is quantifier free.
\item
  Let $K$ be a class of $L$-structures. We call $M \in K$ \emph{existentially closed} for $K$
  if the following holds: for every $M \subseteq N \in K$, a quantifier free formula $\ph(x,\y)$
  and a tuple $\b \in M$, we have $\inf_{x}^M\ph(x,\b) = \inf_{x}^N\ph(x,\b)$.
\item
  If $T$ is a universal theory we say that $M\models T$ is existentially closed for $T$ if
  it is existentially closed for $K = \Mod(T)$. When $T$ is clear from the context we omit it
  and say ``$M$ is existentially closed'' or ``$M$ is an e.c. structure'' or ``$M$ is an e.c.
  model''.
\end{enumerate}
\end{dfn}

\begin{rmk}
  For $M \in K$, to be existentially closed for $K$ means in a sense that anything
  which is quantifier free definable with parameters in $M$, that can happen in some model
  in $K$, happens already in $M$. In this sense, existentially closed models are
  ``generic'' among structures in $K$.
\end{rmk}

\begin{exm}\label{exm:ec}
\begin{enumerate}
\item
  Atomless probability algebras are existentially closed among all probability algebras
  (see \cite{BU} or \cite{BBHU}).
\item
  Hilbert spaces equipped with a unitary operator $U$ with full spectrum
  ($Spec(U) = S^1$) are e.c. among all Hilbert spaces equipped with
  a unitary operator, see \cite{BUZ}.
\item
  Atomless probability algebras with an aperiodic automorphism are e.c. among
  probability algebras equipped with an automorphism, see \cite{BH}.
\end{enumerate}
\end{exm}

Let $T$ be a universal theory, $K = \Mod(T)$, and $K^{ec}$ be the class of e.c. models
of $T$. We call $K^{ec}$ the \emph{continuous Robinson theory} of $T$. One may ask: is
$K^{ec}$ elementary (i.e. is there a continuous theory $T^*$ such that $K^{ec} = \Mod(T^*)$)?
The answer is not always positive, even in the classical (discrete) context. For example,
the Robinson theory of groups (i.e. $T$ is collection of first order sentences which
are true in all groups, $K$ is the class of all groups, and $K^{ec}$ consists of
all groups which are existentially closed) is not elementary.
But often the answer is yes; in this case we say that $T$ admits a model companion and call
$T^*$ the \emph{model companion} of $T$.

It is easy to see that in this case $T^*$ is \emph{model complete}: if $M,N \models T^*$ and
$M\subseteq N$, then $M\prec N$. It does not necessarily
eliminate quantifiers; if it does, we call it the \emph{model completion} of $T$.

\begin{rmk}
  In Example \ref{exm:ec} above, the classes of e.c. models are in fact elementary, and the appropriate
  theories are the model companions, and even the model completions of the universal theories.
\end{rmk}

\begin{obs}
  Let $T$ be a universal theory.
  Then $\fS_T$ is a closed subset of $\fS$, and therefore
  a Polish space.
\end{obs}
\begin{prf}
  Clear.
\end{prf}

The following fact is well-known, but the author is not aware of a
written reference. Although the proof is identical to that of the
classical (discrete) analogue, we include it for completeness. In
order not to scare the reader, we only deal with separable
structures, which is all we need in this article (the proof for an
arbitrary infinite cardinality is essentially the same).

\begin{fct}\label{fct:closedexist}
    Let $T$ be a universal theory, $M\models T$ separable. Then there exists
    a separable $N \supseteq M$, $N \models T$, $N$ is e.c. for
    $T$.
\end{fct}
\begin{prf}
    The proof is standard and resembles very much the construction of the algebraic
    closure of a given field.

    Let $M_0 = M$.
    We construct separable $M_i \models T$
    for $i<\om$ by induction as follows:

    Given $M_i$ let $\seq{\ph_\al(\x_\al) \colon \al<\om}$ be an
    enumeration of all quantifier free formulae $\ph(\x)$ with parameters in $M_i$.
    Now define a sequence $M_i^\al$ of separable models of $T$ by induction on $\al<\om$ as follows:

    \begin{itemize}
    \item
        $M_i^0 = M_i$
    \item
        Given $M_i^\al$, if there is no $\a \in M_i$ satisfying
        $[\ph_\al(\a) = 0]$, but there exists $M' \supseteq M_i^\al$, $M'\models T$
        where such $\a$ exists, let $M_i^{\al+1}$ be any such separable $M'$
        (for the cardinality preservation one can use e.g. Proposition
        7.3 in \cite{BBHU}).
        \\Otherwise let $M_i^{\al+1} = M_i^\al$.
    \end{itemize}
    Now define $M_{i+1} = \overline{\cup_{\al<\lam}M_i^\al}$. Note
    that $M_{i+1}$ has the following property: if there exists a
    quantifier free formula $\ph(\x)$ with parameters in $M_i$, an extension
    $M'\models T$ of $M_{i+1}$ and $\a \in M'$ satisfying $[\ph(\a)
    = 0]$, then such $\a$ exists already in $M_{i+1}$.

    Finally, let $N = \overline{\cup_{i<\om}M_i}$; it is easy to check that it
    is existentially closed.

\end{prf}

\section{Inductive theories}


Recall that we assume that theories are closed under entailment, i.e. every closed condition
which follows from $T$ is already in $T$.
We denote
by $T^o$ the collection of all open conditions which follow from $T$.
Let $T^{oc} = T\cup T^o$.

Let $T$ be an $L$-theory and let $\Delta$ be collection of conditions
(open or/and closed). We denote by $T_\Delta$ the
$\Delta$-\emph{part} of $T$. So  $T_\Delta = T^{oc}\cap\Delta$
For an $L$-structure
$M$, we denote the $\Delta$-part of $\Th(M)$ by $Th_\Delta(M)$.

As usual, we define $\Sigma_n$ and $\Pi_n$ formulae by induction on $n$:
\begin{itemize}
\item
  $\Sigma_0 = \Pi_0 = $ quantifier free formulae
\item
  $\Sigma_{n+1}$ is the collection of formulae of the form $\inf_{\x}\ph(\x,\y)$ where
  $\ph(\x,\y) \in \Pi_n$
\item
  $\Pi_{n+1}$ is the collection of formulae of the form $\sup_{\x}\ph(\x,\y)$ where
  $\ph(\x,\y) \in \Sigma_n$
\end{itemize}

\begin{rmk}
  So $\Sigma_1$ is the collection of all the existential formulae, $\Pi_1$ is the
  collection of all the universal formulae.
\end{rmk}

\begin{dfn}
\begin{enumerate}
\item
  For $\Lambda \subseteq L$, we denote by $\Lambda^o$ the collection of all
  open conditions of the form $\ph<\epss$ for $\ph \in \Lambda$, $\eps>0$.
\item
  For $\Lambda \subseteq L$, we denote by $\Lambda^c$ the collection of all
  closed conditions of the form $\ph\le\epss$ for $\ph \in \Lambda$, $\eps>0$.
\item
  Let $\Delta$ be a collection of conditions. We call a theory $T$ a $\Delta$-\emph{theory}
  if $T_\Delta \models T$.
\item
  We call a theory $T$ \emph{inductive} if it is axiomatizable by open conditions of
  the form $\sup_{\x}\inf_{\y}\ph(\x,\y)<\eps$, where \phh is quantifier-free. So
  $T$ is inductive if it is a $\Pi_2^o$-theory.
\end{enumerate}
\end{dfn}

\begin{rmk}
\begin{enumerate}
\item
  So a theory $T$ is universal iff it is a $\Delta$-theory for $\Delta = \Pi_1^c$.
\item
  Maybe the reader would expect us to work with $\Pi_2^c$-theories instead of $\Pi_2^o$.
  Note that if $T$ is $\Pi_2^c$ then it is inductive, and for complete theories the notions
  are equivalent; but as we want Theorem \ref{thm:ind} to hold for all theories, not
  necessarily complete, and as we want our theories to define $G_\delta$ subsets of
  $\fS$, the natural choice is open conditions.
\end{enumerate}
\end{rmk}




\begin{lem}\label{lem:local}
  Let $\Delta = \Sigma^o_n$ or $\Delta = \Pi^o_n$ for some $n$
  (or just $\Delta$ is a collection of open conditions closed under rescaling, i.e.
  multiplication by scalars and
  the ``pointwise minimum''
  connective).
  Let $T$ be a theory and suppose that
  for every two $L$-structures $M,N$ such that $M \models T$ and
  $\Th_\Delta(M) \subseteq \Th_\Delta(N)$, we have $N \models T$.
  Then $T$ is a $\Delta$-theory.
\end{lem}
\begin{prf}
  Suppose not; so there exists $N \models T_\Delta$, $N \not\models T$.
  By the assumption, for no $M \models T$ do we have $\Th_\Delta(M) \subseteq \Th_\Delta(N)$.
  In other words, for every $M \models T$ there exists a formula $\ph_M$
  with $[\ph_M<\eps] \in \Delta$ such that
  $\ph_M^M<\eps$, $\ph_M^N\ge\eps$. By rescaling we may assume
  $\eps = \frac 12$.

  So the set $T \cup \set{[\ph_M \ge \frac 12]\colon M\models T}$ is inconsistent.
  By compactness, $[\min_{i=1}^k\ph_i<\frac 12] \in T^o$ for some finite
  collection of such $\ph_i$. But $\Delta$ is closed under taking minima, and
  every one of the conditions $[\ph_i<\frac 12]$ is in $\Delta$,
  so (as $N \models T_\Delta$),
  $(\min_{i=1}^k\ph_i)^N<\frac 12$, so for some $i$ $\ph_i^N<\eps$, a contradiction.
\end{prf}


\begin{lem}\label{lem:embed}
  Let $T$ be a complete $L$-theory and $M$ an $L$-structure
  with $\Th_{\Sigma_n^c}(M) \subseteq T_{\Sigma_n^c}$.
  Then there exists $M' \models T$ and a $\Sigma_n$-elementary embedding $f\colon M \to M'$.
\end{lem}
\begin{prf}
  Let $M = \lseq{a}{\al}{\lam}$, $M'$ a $\lam$-saturated model of $T$. Construct
  $f_\al \colon A_\al = \set{a_\be\colon\be<\al} \to M'$ an increasing continuous
  sequence of $\Sigma_n$-elementary embeddings. Given $f_\al$, consider the
  $\Sigma_n$-type in $M$ of $a_\al$ over $A_\al$, call it $\pi(x)$. By the assumption that
  $\Th_{\Sigma_n^c} \subseteq T$, $f(\pi(x))$ is a type in $M'$ over $f(A_\al)$ (as $\Sigma_n$ is closed under
  ``inf''), and use the saturation of $M'$.
\end{prf}

\begin{thm}\label{thm:ind}
  Let $T$ be an $L$-theory such that $\Mod(T)$ is preserved under
  unions of chains, that is if \lseq{M}{i}{\om} is an increasing chain of models of $T$, then
  the closure of the union  $M = \overline{\bigcup_iM_i}$ is also a model of $T$. Then $T$ is inductive.
\end{thm}
\begin{prf}



  Suppose $\Mod(T)$ is closed under unions of chains. We would like to show
  that $T_{\Pi^o_2} \models T$.

  Let $M \models T$ and $T'$ a complete
  $L$-theory extending $T_{\Pi_2^o}$, $N\models T'$ with $\Th_{\Pi_2^o}(M) \subseteq T'$  We will show
  that $N \models T$, which clearly suffices by Lemma \ref{lem:local}.

  Construct a chain
  $N_0 \subseteq M_1 \subseteq N_1 \subseteq \ldots$ such that:
  \begin{enumerate}
  \item
      $N_0 = N$
    \item
      $M_i \models T$, $N_i \models T'$
    \item
      $N_i \prec N_{i+1}$
  \end{enumerate}

  If the construction is possible, we are done:

  Let $M = \overline{\bigcup_iM_i} = \overline{\bigcup_iN_i}$. $M\models T$ by the assumption on $T$ and (ii) above.
  On the other hand, clearly $N_0 \prec M$ (as $N_0 \prec N_i$ for all $i$ by (iii) above), so $N = N_0 \models T$,
  and we are done.

  Why is the construction possible?

  Let $N_0 = N$. As $\Th_{\Pi^o_2}(M) \subseteq \Th(N) = T'$,
  we have $\Th_{\Sigma_2^c}(N) \subseteq \Th(M)$, so by Lemma \ref{lem:embed},
  there exists $M_0 \models \Th(M)$ and a $\Sigma_2^c$-embedding of
  $N_0$ into $M_0$.

  Let $N_0 = \lseq{a}{\al}{\lam}$. Enrich the vocabulary $\tau$ with
  $\lam$-many constant symbols, call the new language $L'$.

  \begin{clm}
  $\Th_{\Sigma_1^c}(M_0,\lseq{a}{\al}{\lam}) \subseteq \Th(N_0,\lseq{a}{\al}{\lam})$
  as $L'$-theories.
  \end{clm}
  \begin{prf}
    Clearly
    (as $N_0$ is a $\Sigma_2^c$-elementary submodel of $M_0$),
    $\Th_{\Sigma_2^c}(N_0,\lseq{a}{\al}{\lam}) \subseteq \Th(M_0,\lseq{a}{\al}{\lam})$
    as $L'$-theories, and therefore
    $\Th_{\Pi_2^o}(M_0,\lseq{a}{\al}{\lam}) \subseteq \Th(N_0,\lseq{a}{\al}{\lam})$
    as $L'$-theories, in particular

    \begin{equation}\label{equ:1}
    \Th_{\Sigma_1^o}(M_0,\lseq{a}{\al}{\lam}) \subseteq \Th(N_0,\lseq{a}{\al}{\lam})
    \end{equation}

    as $L'$-theories.

    Let $[\inf_{\x}\ps(\x,\a)\le\eps]$ be a closed
    existential condition satisfied by $M_0$ with parameters $\a \in N_0$ (i.e.
    $\a = a_{\al_1},\ldots,a_{\al_k}$ for some $\al_1, \ldots, \al_k < \lam$). Then
    $M_0 \models [\inf_{\x}\ps(\x,\a)<\eps']$ for every $\eps' > \eps$. So this is true
    in $N_0$ (by (\ref{equ:1}) above), which completes the proof of the claim.

%
  \end{prf}

  By the Claim above and Lemma \ref{lem:embed}, there exists $N_1 \models T'$ into which
  $M_0$ is $\Sigma_1^c$-embedded in the language $L'$. Clearly, this means that $N_0 \prec N_1$.

  The rest of the construction is similar.

\end{prf}

We obtain the analogue of a well-known Robinson's theorem in the continuous context:

\begin{cor}
  If $T$ is model complete, then it is inductive.
\end{cor}

\begin{cor}
  Let $T$ be a universal theory which has a model companion $T'$. Then $T'$ is
  inductive.
\end{cor}
\begin{prf}
  $T'$ is model complete.
\end{prf}


\section{Generic and random models}
\subsection{Model completions and topological genericity}

\begin{obs}
Let $T'$ be an inductive theory. Then $\fS_{T'}$ is a $G_\delta$ subset of $\fS$.
\end{obs}
\begin{prf}
  For every quantifier free formula $\ph(\x)$, $\a,\b \in \setN$ and $\eps>0$, the open condition
  $[\ph(\a,\b) < \eps]$ defines an open subset of $\fS$, which we called $U_{\ph(\a,\b),\eps}$.
  The open condition $[\inf_{\x}\ph(\x,\b) < \eps]$
  corresponds, therefore, to an open subset of $\fS$, which equals
  $\bigcup_{\a\in \setN}U_{\ph(\a,\b),\eps}$. A $\Pi^o_2$ condition defines a subset of
  $\fS$ which is a (countable) intersection (over all possible $\b \in \setN$) of
  sets as above; therefore it is a $G_\delta$ set. Clearly, a countable collection of
  $\Pi^o_2$ conditions still corresponds to a $G_\delta$ set.
\end{prf}

\begin{fct}
  Let $T$ be a universal theory. Then the collection of separable e.c. models is \emph{dense} in
  $\fS_T$.
\end{fct}
\begin{prf}
  By Fact \ref{fct:closedexist} every separable $M\models T$ can be extended to a separable
   e.c. model $M'\models T$. Now it is easy to see that one can rename the elements of $M'$
  such that a certain finite $\a \in \setN$ remains unchanged (and so $M'$ is indeed in a specified
  open neighborhood of $M$ in $\fS_T$).
\end{prf}

\begin{cor}\label{cor:gdeltadense}
  Let $T$ be a universal theory which has a model companion $T'$. Then $\fS_{T'}$ is
  a $G_\delta$ dense subset of $\fS_T$.
\end{cor}

Recall that a theory $T$ is called $\aleph_0$-\emph{categorical} if any two separable
models of $T$ are isomorphic.

\begin{dfn}\label{dfn:generic}
Let $T$ be a universal theory. We call $M \in \fS_T$ \emph{generic} if
the isomorphism class of $M$ is $G_\delta$
dense in $\fS_T$.
\end{dfn}

\begin{cor}\label{cor:generic}
  Let $T$ be a universal theory which has a model companion $T'$, and assume $T'$ is
  $\aleph_0$-categorical. Then (any) existentially closed model of $T$ is generic in
  $\fS_T$.
\end{cor}
\begin{prf}
  By Corollary \ref{cor:gdeltadense}, $\fS_{T'}$ is $G_\delta$ dense in $\fS_T$. By $\aleph_0$-
  categoricity of $T'$, $\fS_{T'}$ is the isomorphism class of any e.c. model of $T$ (which is in $\fS_T$).
\end{prf}


\subsection{Random structures}

Once we have shown that the class of existentially closed models in $\fS$ is ``big'' in the sense of Baire
category theory, a natural question is: is there a similar measure-theoretic result?
In \cite{Ver02}
Vershik shows that the Urysohn space is in a sense the random metric space. We know that the model
companion of the universal theory of graphs is the random graph. Are these facts  particular cases
of a model theoretic phenomenon?

Let $T$ be a universal theory, $\mu$ a probability measure on $\fS_T$ satisfying

\begin{asm}\label{asm:measure}
\begin{enumerate}
\item No nonempty open set has probability $0$.
\item $\mu$ is invariant under the action of $S_\infty$ on $\fS_T$. In other words,
for every formula $\ph(\x)$, $\eps>0$ and $\a,\b \in \setN$, we have the equality
$\mu(U_{\ph(\a),\eps}) = \mu(U_{\ph(\b),\eps})$. So $\mu(U_{\ph(\a),\eps})$ does not depend
on $\a$.
\end{enumerate}
\end{asm}

Clearly, these are very natural assumptions on a measure on $\fS_T$, once we are interested
in ``random structures'': first, we assume that if a certain open event occurs in some model of $T$,
then its probability is positive. Second, we assume that in a sense isomorphic models ``occur''
with equal probability.

\begin{lem}
  Let $\mu$ be as above.
  Then the set of all existentially closed structures in $\fS_T$ has probability $1$. In other
  words, if we pick a structure ``randomly'', it is going to be
  existentially closed almost surely.
\end{lem}
\begin{prf}
  Let $M$ be a randomly chosen structure. We aim to show that
  with
  probability $1$ it is existentially closed.
  Let $\ph(\x)$ be a formula, and suppose that in some $M\subseteq N \models T$
  we have $\inf_{\x}^N\ph(\x) \le \eps$. Let $\eps'>\eps$. So there exists
  $\a \in N$ such that $\ph^N(\a)<\eps'$, and therefore $\mu(U_{\ph(\a),\eps'}) = \delta > 0$.
  By the invariance of $\mu$, $\mu(U_{\ph(\b),\eps'}) = \delta$ for every $\b \in \setN$, and
  so the probability that in a randomly chosen structure $M$ we have
  $\ph^M(\b)\ge \eps'$ is bounded away from $1$ for each $\b \in M$.

  Now clearly with probability $1$ for some $\b \in M$ we have $\ph^M(\b)<\eps'$, therefore
  $\inf_{\x}^M\ph(\x) \le \eps$ almost surely, and we are done.
\end{prf}

\begin{rmk}
  Note that we did not really use the invariance of $\mu$. We only need that the probability
  of the event $U_{\ph(\a),\eps}$ is either $0$ for all $\a$ or bounded away from $0$ for all
    $\a$.
\end{rmk}

\begin{cor}\label{cor:ecrandom}
  Let $T$ be a universal theory which has a model companion $T'$. Then $\fS_{T'}$ is a set
  of probability $1$ in $\fS_T$.
\end{cor}

\begin{dfn}\label{dfn:random}
  We call a separable model of a universal theory $T$ \emph{random} if the measure of its
  isomorphism class in $\fS_T$ is $1$ with respect to any probability measure $\mu$ as in
  Assumption \ref{asm:measure}. In other words, $M$ is a random model of $T$
  if for every $\mu$ as above, a randomly chosen
  structure in $\fS_T$ is almost surely isomorphic to $M$.
\end{dfn}

Just like in Corollary \ref{cor:generic} we obtain:

\begin{cor}\label{cor:random}
  Let $T$ be a universal theory which has a model companion $T'$.
  Assume furthermore that $T'$ is $\aleph_0$-categorical. Then any separable
  model of $T'$ is a random model of $T$.
\end{cor}

Clearly, this generalizes the ``randomness'' of the countable random graph; see more
in the following subsection.


\subsection{Concluding remarks on genericity}

In this section we have shown that the model theoretic notion of genericity gives rise
to both Baire category theoretical and measure theoretical notions of genericity in the space $\fS$.
In other words, we have shown:

\begin{cor}
Let $T$ be a universal theory which admits a model companion $T'$.
Then $\fS_{T'}$ is both $G_\delta$ dense in $\fS$ and of measure $1$ with respect to any
reasonable measure on $\fS$ (i.e. any measure satisfying Assumption \ref{asm:measure}).
\end{cor}

In particular, we have the following:

\begin{cor}\label{cor:genericrandom}
Suppose $T$ is a universal theory which has a model companion $T'$,
and assume furthermore that $T'$ is $\aleph_0$-categorical. Then the
(unique up to isomorphism) model of $T'$ is both the generic and the
random model of $T$.
\end{cor}

\begin{exm}
  The atomless separable probability algebra is both the generic and
  the random separable probability algebra.
\end{exm}
\begin{prf}
  The theory of atomless probability algebras is the model companion of the universal
  theory of probability algebras by \cite{BU}. It is also $\aleph_0$-categorical, so
  apply Corollary \ref{cor:genericrandom}.
\end{prf}

As we have already mentioned, every classical first order theory is a continuous first
order theory with discrete metric. We can therefore apply our analysis to e.g. the
theory of the random graph. Recall that the theory of the random graph is the model
completion of the universal theory of graphs.

\begin{exm}
  The random graph is the generic countable graph.
\end{exm}
\begin{prf}
  The (classical) first order theory of the random graph is the model completion (and therefore
  the model companion) of the universal theory of graphs. It us also $\aleph_0$-categorical.
  So the desired conclusion follows from Corollary \ref{cor:genericrandom}.
\end{prf}

Similarly, the unique countable model of the model completion of the universal
theory of graphs is the random graph in the sense defined here in Definition
\ref{dfn:random}.
Well, no surprise here: we're just saying that the random graph is, well, random.

In the following section we will show that the continuous first order theory of the Urysohn space
has similar properties, and therefore Corollary \ref{cor:genericrandom} applies to it as well. One can
think of this theory as the continuous analogue of the theory of the random graph: instead of
the discrete predicate $R(x,y)$ in the theory of graphs which can be either true or false, we have
a metric which can take any value between $0$ and $1$.

\section{Urysohn space}

Many results on the model theory of the Urysohn space here are ``folklore'', but the author
is not aware of any written references. In order to follow the proofs, the reader
should be familiar with basics of continuous
model theory slightly beyond what is sketched in section 2 of the article.

We remind the reader that the \emph{Urysohn space} is the universal complete separable metric space, first
constructed by Pavel Urysohn. Due to the limitations of the genre, we will consider the \emph{bounded}
Urysohn space, i.e. Urysohn space of diameter 1. We denote it by $\U$.

Denote by $\frak E_n$ the collection of all possible distance
configurations on $n$ points of diameter 1. It will be convenient
for us to think about it in the following way: $\vartheta(x_1,
\ldots, x_n) \in \frak E_n$ if $\vartheta(x_1, \dots, x_n)$ is a
formula of the form $$\bigvee_{1\le i,j\le
k}\left|d(x_i,x_j)-r_{ij}\right|$$ where the matrix $(r_{ij})_{1\le
i,j\le k}$ is a distance matrix of some finite metric space of
diameter 1, and $\bigvee$ stands for the lattice operation of
pointwise maximum.

Let us introduce the following notation: for $\vartheta \in \frak E_{n+1}$, let $\vartheta \rest n$
be the restriction of $\vartheta$ to the first $n$ variables.

Clearly, for every $\vartheta \in \frak E_{n+1}$, for every $\eps>0$ there exists a $\delta = \delta(\eps) >0$
such that if $a_1, \ldots, a_n \in \U$ satisfy $\vartheta \rest n(a_1, \ldots, a_n) < \delta$, then
there exists $a_{n+1} \in \U$ such that $\vartheta(a_1, \ldots, a_n, a_{n+1}) \le \eps$.

Let $T_\U$ be the collection of all the conditions of the form
\begin{equation*}
  \left[\sup_{x_1, \ldots, x_n}\inf_{y}
    \left( \frac{\eps}{1-\delta} \left(1-\vartheta\rest n(x_1, \ldots, x_n) \right) \bigwedge
    \vartheta(x_1, \ldots, x_n,y) \right) \le \eps \right]
\end{equation*}

which is just one way of stating

\begin{equation*}
  \forall x_1, \ldots x_n \exists y (\vartheta\rest n(x_1, \ldots, x_n)<\delta \rightarrow
  \vartheta(x_1,\ldots,x_n,y)\le \eps)
\end{equation*}

Note that $\bigwedge$ stands for the lattice operation of pointwise minimum.

The following follows from the standard Urysohn's argument:

\begin{fct}
  The only separable complete model of $T_\U$ is $\U$.
\end{fct}

\begin{cor}
  $T_\U$ is $\aleph_0$-categorical, and therefore a complete continuous theory.
\end{cor}
\begin{prf}
  By (the continuous version of) Vaught's test.
\end{prf}

\begin{prp}
  $T_\U$ eliminates quantifiers.
\end{prp}
\begin{prf}
  By the classical back-and-forth argument (see Theorem 4.16 in \cite{BU}) using the
  axioms of $T_\U$.
\end{prf}

\begin{cor}\label{cor:urysohn}
  $T_\U$ is the model completion (and therefore the model companion) of
  the ``empty'' continuous universal theory (the universal theory of
  a metric space with no extra-structure). $\U$ is (the only) existentially
  closed metric space.
\end{cor}

A natural conclusion from our analysis is the following form of Vershik's theorems:

\begin{cor}
  The Urysohn space (of diameter 1) is the generic and the random metric space (of diameter 1).
\end{cor}
\begin{prf}
  The theory of the Urysohn space is the model companion of the universal
  theory of metric spaces and is $\aleph_0$-categorical, so the result
  follows immediately from Corollary \ref{cor:genericrandom}.
\end{prf}


\end{document}